\documentclass[12pt]{article}   
\usepackage{amsmath}
\usepackage{amssymb}
\usepackage{amsthm}
\usepackage[totalwidth=17cm, totalheight=25cm]{geometry}
\usepackage{graphicx}
\usepackage{mathabx}
\usepackage{tabularx}
	\newcolumntype{C}[1]{>{\centering\arraybackslash}m{#1}} 
	\newcolumntype{R}[1]{>{\raggedleft\arraybackslash}m{#1}} 

\newtheoremstyle{boldplain}
{9pt}
{9pt}
{\itshape}
{}
{\bfseries}
{.}
{.5em}
{\thmname{#1}\thmnumber{ #2}\thmnote{ (#3)}}%

\newtheoremstyle{bolddefinition}
{9pt}
{9pt}
{}
{}
{\bfseries}
{.}
{.5em}
{\thmname{#1}\thmnumber{ #2}\thmnote{ (#3)}}%

\theoremstyle{boldplain}

\newtheorem{cor}[equation]{Corollary}

\newtheorem{lem}[equation]{Lemma}
\newtheorem{lemma}[equation]{Lemma}
\newtheorem{prop}[equation]{Proposition}

\newtheorem{question}[equation]{Question}

\newtheorem{thm}[equation]{Theorem}

\theoremstyle{bolddefinition}

\newtheorem{rem}[equation]{Remark}

\newfont{\bigbf}{cmbx10 scaled\magstep1}
\normalbaselines
\def\R{{\mathbb R}}
\def\H{{\mathbb H}}

\def\N{{\mathbb N}}
\def\P{{\mathbb P}}

\def\ga{\gamma}
\def\Ga{\Gamma}

\def\la{\lambda}
\def\La{\Lambda}
\def\si{\sigma}

\def\Om{\Omega}

\def\3{\ss}
\def\Dt{\partial_{\tau_{mod}}}

\def\acts{\curvearrowright}
\def\amod{a_{mod}}

\def\Flagt{\operatorname{Flag_{\tau_{mod}}}}

\def\Fmod{F_{mod}}

\def\geo{\partial_{\infty}}

\def\half{\frac{1}{2}}

\def\id{\mathop{\hbox{id}}}

\def\Lat{\La_{\tau_{mod}}}

\def\ol{\overline}

\def\pihalf{\frac{\pi}{2}}

\def\2pithird{\frac{2\pi}{3}}

\def\pos{\mathop{\hbox{pos}}\nolimits}

\def\simod{\si_{mod}}

\def\span{\operatorname{span}}

\def\st{\operatorname{st}}

\def\ost{\operatorname{ost}}

\def\tangle{\angle_{Tits}}
\def\taumod{\tau_{mod}}

\def\tits{\partial_{Tits}}
\def\Th{\mathop{\hbox{Th}}\nolimits}

\def\8{\infty}
\def\<{\langle}
\def\>{\rangle}

\def\BI{\begin{itemize}}
\def\EI{\end{itemize}}

\long\def\comment#1\endcomment{}

\def\F{\mathrm F}

\title{Domains of discontinuity of Lorentzian affine group actions}
\author{Michael Kapovich and Bernhard Leeb}
\date{\today}							

\begin{document}
\maketitle


\begin{abstract}
We prove nonemptyness of domains of proper discontinuity  of Anosov groups  
of affine Lorentzian transformations of $\R^n$. 
\end{abstract}

\section{Introduction}

There is a substantial body of literature, going back to the pioneering work of Margulis \cite{Margulis},  
on properly discontinuous non-amenable groups of affine transformations, see e.g. \cite{Abels, AMS02, AMS11, Drumm, DGK, GLM, Mess}, and numerous other papers, in particular, the recent survey \cite{DDGS}. In this paper we address a somewhat related question of nonemptyness of domains of proper discontinuity of discrete groups acting on  affine spaces:

\begin{question}
Which discrete subgroups $\Ga< Aff(\R^n)$ have nonempty discontinuity domain in the affine space $\R^n$? 
\end{question}

\noindent In this paper we limit ourselves to the following setting: Suppose that $\Ga< \R^n \rtimes O(n-1,1)< Aff(\R^n)$ 
is a discrete subgroup such that the linear projection $\ell: \Ga\to O(n-1,1)$ is 
a {\em faithful representation with convex-cocompact image},  
see e.g. \cite{Bowditch} for the precise definition. Given a representation $\ell: \Ga\to O(n-1,1)$, the  affine 
action of $\Ga$ is determined by a cocycle $c\in Z^1(\Ga, \R^{n-1,1}_\ell)$.  
Even in the case $n=3$ and $\ell(\Ga)$ a Schottky subgroup of $O(2,1)$ (which is the setting of Margulis' original examples), while some actions are properly discontinuous on the entire $\R^3$ (as proven by Margulis, see also \cite{GLM} for  a general description of such actions), nonemptyness of domains of discontinuity for {\em arbitrary} $c$ does not appear to be obvious\footnote{The reaction to the question that we observed included: ``clearly true", ``clearly false", ``unclear".}.

  The main result of this note is:

\begin{thm}\label{thm:main}
Every subgroup $\Ga< \R^n \rtimes O(n-1,1)$ with faithful convex-cocompact  linear representation 
$\ell: \Ga\to O(n-1,1)$, acts properly discontinuously on a nonempty open subset of the Lorentzian space $\R^{n-1,1}$. 
\end{thm}

We will prove this theorem by applying results on domains of discontinuity for discrete group actions on flag-manifolds proven in \cite{coco15}. More precisely, we will check that $\Ga$ is a {\em $\taumod$-Anosov subgroup} of the Lie group $G=O(n,2)$ for a suitable model simplex (actually, a vertex) $\taumod\subset \simod$. In section 
\ref{sec:identification} we will equivariantly identify the Lorentzian space $\R^{n-1,1}$  and an open Schubert cell in a partial flag-manifold 
$\F_1=G/P_{\taumod}$ of the group $G=O(n,2)$. In \cite{coco15} we proved that for each $\taumod$-Anosov subgroup $\Gamma$ of a semisimple Lie group $G$ and each {\em fat thickening} $\Th(\Lat(\Ga))\subset \F_1$ of the $\taumod$-limit set $\Lat(\Ga)\subset \F_1$, the group $\Ga$ acts properly discontinuously on the open subset $\Omega_{{\scriptsize \Th}}(\Ga)= \F_1\setminus  \Th(\Lat(\Ga))$. In Section \ref{sec:main} of this paper we verify that $\Omega_{{\scriptsize \Th}}(\Ga)\ne \emptyset$ in the context of  $\taumod$-Anosov subgroups 
$\Ga< \R^n \rtimes O(n-1,1)< O(n,2)$ and the {\em maximal thickening} $\Th$. This, in turn, will establish nonemptyness of the domain of discontinuity of $\Gamma$ in $\R^{n-1,1}$.

\section{Geometric preliminaries}\label{sec:prelim}

{\bf Symmetric spaces of noncompact type and their visual boundaries.} For basics of symmetric spaces and their visual boundaries we refer the reader to \cite{BGS} and \cite{Eberlein}. 

Consider a symmetric space of noncompact type $X=G/K$, where $G$ is a semisimple Lie group (with finite center) and $K$ is its maximal compact subgroup. Fix also a base-point $o\in X$ (the choice is ultimately irrelevant), fixed by $K$. We let $d$ denote the Riemannian distance function on $X$ and $\angle_x(y,z)$ the Riemannian angle between nondegenerate geodesic segments $xy, xz$ emanating from $x$. The {\em visual boundary} $\geo X$ of $X$, as a set, is identified with the set of equivalence classes $[\rho]$ of geodesic rays $\rho: \R_+\to X$ in $X$, where two rays are equivalent if and only if their images are at a finite Hausdorff distance from each other. One says that every ray $\rho$ representing $\xi=[\rho]$ is {\em asymptotic to $\rho$.} The {\em Tits angle} $\tangle(\xi_1, \xi_2)$ between points $\xi_1=[\rho_1], \xi_2=[\rho_2]$ is defined as 
$$
\sup_{x\in X}  \angle_x(\rho_1(t), \rho_2(t)),
$$ 
 where the supremum is taken over all pairs of rays $\rho_1, \rho_2$ representing $\xi_1, \xi_2$ such that $\rho_1(0)=\rho_2(0)=x$. 
 Since $X$ is a symmetric space, there exists a flat $F\subset X$ such that $\xi_1, \xi_2$ are represented by rays whose images are contained in $F$. The supremum in the definition of $\tangle(\xi_1, \xi_2)$ is realized by   pairs of such rays.  The Tits angle defines the {\em Tits metric} on $\geo X$. This metric is invariant under the natural $G$-action on $\geo X$. 
 
 The visual boundary of $X$ has two natural topologies. The first one is the {\em visual topology}: Every $\xi\in \geo X$ is represented by a unique unit speed geodesic ray emanating from $o$. Thus, there is a natural bijection between $\geo X$ and the unit sphere in the tangent space $T_oX$. The visual topology on $\geo X$ is the one making this bijection a homeomorphism. The natural $G$-action on $\geo X$ is continuous with respect to this topology. This topology extends to a {\em visual compactification} $\overline{X}=X\cup \geo X$: A sequence $(x_n)$ in $X$ converges to $\xi=[\rho]\in \geo X$ if 
 $$
 \lim_{n\to\infty} \angle_o(x_n, \rho(1))=0 \quad \hbox{and}\quad \lim_{n\to\infty} d(o,x_n)=\infty, 
 $$ 
 where $\rho(0)=o$. For a subset $A\subset X$, the {\em visual boundary} of $A$ is the intersection of $\geo X$ with the closure of $A$ in 
 $\overline{X}$ with respect to the visual topology. 
 
  The second, {\em Tits topology}, is the one defined by the Tits metric. With respect to this topology,  $\geo X$ has the structure of a certain simplicial complex, the {\em spherical (Tits) building} $\tits X$, invariant under the action of $G$. We will fix a {\em model chamber} of $\geo X$, i.e. a  facet $\simod$ of this spherical  building, a {\em model maximal flat} $\Fmod\subset X$, it is the unique maximal flat in $X$ whose visual boundary $\amod$ (the {\em model apartment} in $\tits X$) is a subcomplex containing  $\simod$ and such that $o\in \Fmod$.  
 The {\em Euclidean Weyl chamber} $\Delta$ of $X$ is the cone in $\Fmod$ with the tip $o$ over  
 $\simod$ (the union of geodesic rays emanating from $o$ and asymptotic to the points of $\simod$). The {\em Weyl group} $W$ of $X$ 
 is the image of $K\cap Stab_G(\Fmod)$ in the isometry group of the flat $\Fmod$. Then $\Delta$ is a fundamental domain of the $W$-action on $\Fmod$.  The Weyl group $W$ has a standard word-metric; we let $w_0$ denote the unique longest element of $W$ with respect to this metric. Identifying $\Fmod$ with $\R^r$ (where $r$ is the rank of $X$), we get the {\em opposition involution} $\iota=-w_0$ preserving $\simod$.  In the case of symmetric spaces of type $B$, as in this paper, $w_0=-\id$ and, accordingly, $\iota=\id$. 

{\bf Antipodality.} Two points $\xi, \eta$ in $\geo X$ are called {\em opposite} if $\tangle(\xi, \eta)=\pi$, equivalently, if there exists a geodesic $c$ in $X$ whose opposite subrays are asymptotic to $\xi$ and $\eta$ respectively. Equivalently, there exists a Cartan involution of $X$ swapping $\xi$ and $\eta$. Two simplices $\tau, \hat\tau$ in $\tits X$ are {\em opposite} (or, {\em antipodal}) if and only if they contain opposite {\em generic} points in $\geo X$. (A point in a simplex $\tau$ is {\em generic} if it does not belong to any proper face of $\tau$.) Two simplices in $\tits X$ are opposite  if and only if they are swapped by a Cartan involution of $X$. 

\medskip
{\bf Horoballs.} For every point $\xi=[\rho]$ in $\geo X$ one defines the {\em Busemann function} $b_\xi$ on $X$ (or, more precisely, a family of Busemann functions which differ by additive constants):
$$
b_\xi(x)= \lim_{t\to\infty} (d(\rho(0), x) - t).
$$
Busemann functions satisfy the following equivariance condition with respect to the action of isometries $g$ of $X$:
$$
b_{g\xi}= b_\xi \circ g + Const. 
$$
Sublevel sets of Busemann functions $b_\xi$ are called {\em horoballs centered at $\xi$} and denoted $Hbo$. Busemann functions and, hence,  horoballs, are convex. We will need the following lemma that can be found in \cite[Lemma 4.10]{BGS} and 
 \cite[Proposition 3.4.3]{Eberlein}: 

\begin{lemma}\label{lem:horoball} 
For each horoball $Hbo$ in $X$ centered at $\xi$, the visual boundary of $Hbo$ equals the closed $\pihalf$-ball $\bar{B}(\xi, \pihalf)$ 
in $\geo X$ centered at $\xi$, where the distance is computed in the Tits metric on $\geo X$. 
\end{lemma}

\medskip
{\bf Parallel sets.} Fix two opposite points $\xi, \hat\xi\in \geo X$. The {\em parallel set} $P(\xi, \hat \xi)$ is a certain symmetric subspace in $X$, which is the union of all geodesics $l$ in $X$ that are forward-asymptotic to $\xi\in \geo X$ and 
backward-asymptotic to $\hat{\xi}\in \geo X$. Suppose that $\xi, \hat\xi$ are generic points of two opposite simplices $\tau, \hat\tau$ in $\tits X$. Then $P(\xi, \hat \xi)$ splits isometrically as a direct product $F_{\tau,\hat\tau}\times Y$, where $F_{\tau,\hat\tau}$ is a flat in 
$X$ of dimension $\dim(\tau)+1$ and $Y$ is a (totally-geodesic) symmetric subspace of noncompact type in $X$, called {\em a cross-section} of $P(\xi, \hat \xi)$. In the case of interest to us, $\tau, \hat\tau$ are vertices in $\tits X$, $F_{\tau,\hat\tau}$ is 1-dimensional and $Y$ is a symmetric space of rank 1 (actually, the hyperbolic space). The pointwise stabilizer $G_{\tau,\hat\tau}$ of $\{\tau, \hat\tau\}$ is a reductive subgroup of $P_\tau$; it splits off as a product $G_Y \times \R^r$, where  $\R^r$ is the group of transvections in $G$ preserving the flat $F_{\tau,\hat\tau}$ and 
$G_Y$ is a semisimple Lie group, it is the stabilizer of $Y$ in $G_{\tau,\hat\tau}$. The action of $G_Y$ on $Y$ (and $P(\xi, \hat \xi)$) 
may have a nontrivial (but compact) kernel and the image of $G_Y$ in the isometry group of $Y$ is a subgroup of finite index. The {\em unipotent radical} $U_\tau\triangleleft P_{\tau}$  is a normal subgroup such that $P_\tau=  U_\tau \rtimes G_{\tau,\hat\tau}$. A  more refined form of this decomposition is
$$
P_\tau=  (U_\tau \rtimes G_Y) \rtimes \R^r.
$$
The subgroup $U_\tau \rtimes G_Y$ preserves each horoball centered at $\xi$. See \cite[\S 2.8, 2.10]{anolec} for more details.

\medskip 
For the material below we refer the reader to \cite{coco15} and \cite{manicures}.   

For each point $x\in X$ one defines the $\Delta$-valued distance $d_\Delta(o,x)$ as the unique point of intersection $Kx\cap \Delta$. 
(This definition extends to general pairs of points in $X$ by $G$-invaraince.) 
Consider a face $\taumod$ of the spherical Weyl chamber $\simod$ of $X$. These faces parameterize {\em standard parabolic subgroups} $P_{\taumod}$ of $G$, their $G$-stabilizers. 
The {\em $\taumod$-boundary} $\Dt\simod$ of 
$\simod$ is the union of the faces of $\simod$ which do not contain $\taumod$. The {\em open star} $\ost(\taumod)$ 
of $\taumod$ in $\simod$ is the complement $\simod \setminus \Dt\simod$. 
In the example relevant to us, when $\simod=[u,v]$ is a simplex with the vertices $u, v$ and $\taumod$ is one of the vertices of $\simod$, say, $u$, 
$\Dt\simod=\{v\}$ and $\ost(\taumod)= [u,v)= \simod \setminus \{v\}$. In general, one defines $V(\Dt\simod)\subset \Delta$ as the cone over $\Dt\simod$. 

\medskip {\bf Stars at infinity.} The group $G$ acts transitively on the set of facets of $\geo X$; thus, a face $\tau$ of $\geo X$ is said to have the type $\taumod$ if they lie in the same $G$-orbit. One defines {\em open stars} $\ost(\tau)$ of faces $\tau$ of $\geo X$: One first takes its {\em star},  $\st(\tau)$, the subcomplex in $\tits X$ which is the union of faces containing $\tau$, and then removes from  
$\st(\tau)$ those faces which do not contain $\tau$. In the case of interest to us, $\geo X$ is 1-dimensional (a connected 
graph of valence continuum at each vertex), $\tau$ is a vertex of $\geo X$, $\st(\tau)$  is the union of edges (including their respective vertices!) containing $\tau$ as an end-point and $\ost(\tau)$ is the interior of $\st(\tau)$ with respect to the Tits topology, i.e. the topology of the graph $\tits X$.  A point $\xi\in \geo X$ is said to be {\em $\taumod$-regular} if it belongs to $\ost(\tau)$ for some 
$\tau\in \Flagt$. One quantifies this notion of regularity by taking a compact subset $\Theta\subset \ost(\taumod)$;  a 
 $\xi\in \geo X$ is said to be $\Theta$-regular if its projection to $\simod$ belongs to $\Theta$. 

\medskip
{\bf Flag-manifolds.} Fix a model simplex $\taumod$.  The $G$-orbit  $G \taumod$ is naturally identified with the quotient 
$G/P_{\taumod}$. From the viewpoint of the Tits topology, this quotient is discrete, but, it also has a natural manifold topology (the quotient topology of the Lie group $G$), making it a {\em partial flag-manifold} $\Flagt$. Another way to describe this topology is to note that 
there is a $G$-equivariant bijection between  $G/P_{\taumod}$ and the orbit $G\xi$ for a generic point $\xi\in \taumod$. 
This bijection is a homeomorphism from  $G/P_{\taumod}$ to $G\xi$, where the latter is equipped with the  
subspace topology inherited from the visual topology on $\geo X$.

 \medskip
{\bf Thickenings.} 
We fix a model face $\taumod$ of $\simod$. The $W$-orbit of $\taumod$ in the model apartment 
$\amod$ is naturally identified with the quotient 
$W/W_{\taumod}$, where $W_{\taumod}$ is the stabilizer of $\taumod$ in $W$. The group $W$ acts on $W/W_{\taumod}$ via the left multiplication. The {\em strong Bruhat order} $\le$ on $W$ descends to the {\em folding (partial) order} on $W/W_{\taumod}$: 

$[w]\le [w']$ if and only if representatives $w, w'$ or $[w], [w']\in  W/W_{\taumod}$ can be chosen so that $w\le w'$. 

\medskip 
An {\em ideal} in the poset $(W/W_{\taumod}, \le)$ is a proper subset (i.e., a nonempty subset with nonempty complement) $I$ satisfying the property that with every $[w]\in I$, the ideal contains all smaller elements of $W/W_{\taumod}$. The poset $(W/W_{\taumod}, \le)$ has a unique maximal element $[w_0]$ where $w_0$ is the longest element of $W$. Accordingly, $(W/W_{\taumod}, \le)$ has a unique  
maximal ideal $J$ equal to the complement of $\{[w_0]\}$.  An ideal $I$ is called {\em fat} if 
$$
I \cup w_0 I = W/W_{\taumod}. 
$$
For instance, the unique maximal ideal is fat. 

For every pair of simplices $\tau, \tau'\in\Flagt$, there exist $g\in G$ and $w\in W$ such that $g(\tau)=\taumod$ and $g(\tau')= \nu=w\taumod$, a simplex in $\amod$. The simplex $\nu$ is not uniquely determined by this, but its $W_{\tau}$-orbit is uniquely determined. Hence, we define the {\em relative position} of $\tau'$ with respect to $\tau$, $\pos(\tau', \tau)$, as the 
$W_{\tau}$-orbit of $\nu$, equivalently, the corresponding $W_{\taumod}$-orbit in $W/W_{\taumod}$ (or, equivalently, 
the double coset of $w$ in $W_{\taumod}\backslash W/W_{\taumod}$). 
Let  $I\subset W/W_{\taumod}$ be an ideal invariant under the left $W_{\taumod}$-action. (For instance, the unique maximal ideal satisfies this condition.) For a simplex $\tau\in \Flagt$, we define the {\em thickening} $\Th(\tau)=\Th_I(\tau)\subset \Flagt$ as the subset consisting of simplices $\tau'$ such that $\pos(\tau',\tau)\subset I$. In other words, $\tau'\in \Th(\tau)$ if and only if there exists $g\in G$ such that $g(\tau)=\taumod$ and $g(\tau')\in I$. The thickening $\Th(\tau)$ is a certain closed subcomplex (a union of Schubert cycles)  in a cellular decomposition of $\Flagt$ relative to $\tau$. The thickenings $\Th(\tau)$ satisfy 
$$
\Th(g \tau)= g \Th(\tau), g\in G. 
$$
Given a subset $A\subset \Flagt$ and a $W_{\taumod}$-invariant ideal $I$ in $(W/W_{\taumod}, \le)$, we define the corresponding thickening of $A$ as 
$$
\Th(A)=\bigcup_{\tau\in A} \Th(\tau). 
$$
It is observed in \cite{coco15} (see also \cite[Lemma 8.18]{bordif} and Lemma \ref{lem:L} of this paper) that for every closed subset 
$A\subset \Flagt$ and an ideal $I$, the corresponding thickening $\Th(A)$ is a closed subset of $\Flagt$.   A thickening is called {\em fat} if the corresponding ideal in  $W/W_{\taumod}$ is fat. A thickening is {\em maximal} if the corresponding ideal is the maximal ideal. 

\medskip
{\bf Regularity and flag-convergence.} A nondegenerate geodesic segment $xy$ in $X$ is said to be {\em $\taumod$-regular} if 
$d_\Delta(x,y)\in \ost(\taumod)$. 

A sequence $(x_n)$ in $X$ is said to be {\em $\taumod$-regular} if the sequence of vectors $d_\Delta(o, x_n)\in \Delta$ diverges away from $V(\Dt\simod)$ as $n\to\infty$. In the example relevant to us, when $G$ has rank two and, accordingly, $\Delta$ is two-dimensional, and $\taumod$ is a vertex of an edge $\simod$,  $V(\Dt\simod)$ is the null-set of a certain linear functional on $\Delta$, a simple root $\alpha$. Then $\taumod$-regularity of $(x_n)$ means that 
$$
\lim_{n\to\infty} \alpha \left(d_\Delta(o, x_n)\right)= \infty. 
$$

A sequence $(x_n)$  is said to be {\em uniformly $\taumod$-regular} if  the sequence of vectors $d_\Delta(o, x_n)\in \Delta$ diverges away from $V(\Dt\simod)$ at a linear speed with respect to $d(o,x_n)$. In a more quantitative way, one describes uniformly regular sequences as follows. Fix a compact subset $\Theta\subset \ost(\taumod)$. A  sequence $(x_n)$  is said to be 
{\em $\Theta$-regular} if $d(o,x_n)\to\infty$ and for all but finitely many members of the sequence, the geodesic rays $\rho_n$ from $o$ through $d_\Delta(o, x_n)$ are asymptotic to points of $\Theta$. Then a sequence $(x_n)$ is  uniformly $\taumod$-regular if and only if it is $\Theta$-regular for some compact $\Theta\subset \ost(\taumod)$. 

The same definitions apply to sequences $(g_n)$ in $G$: A sequence $(g_n)$ is (uniformly) $\taumod$-regular if for some (equivalently, every) $x\in X$, the sequence $x_n=g_n(x)$ is (uniformly) $\taumod$-regular. 

In \cite{mlem} we defined a partial compactification of $X$, $\ol{X}^{\taumod}= X\cup \Flagt$. Below we will only describe the notion of {\em flag-convergence} for $\taumod$-regular sequences in $X$ to points of $\Flagt$ with respect to the topology of  
$\ol{X}^{\taumod}$. If $X$ has rank 1, then $\simod$ is a singleton, $\taumod=\simod$  and $\Flagt=\geo X$ (with the visual topology). Accordingly, a sequence $(x_n)$ converges to $\tau\in \Flagt$ if and only if it converges to $\tau\in \geo X$ in the visual topology. 

In higher rank, a ray geodesic $o\xi_n$ through $x_n$ need not even terminate in a face $\tau_n$ of $\tits X$ of type $\taumod$. But, if it does,  then $x_n\to \tau\in \Flagt$ if and only if $\tau_n\to \tau$ in $\Flagt$. 

In general, one defines flag-convergence $x_n\to \tau\in \Flagt$ for $\taumod$-regular sequences $(x_n)$ in $X$ as follows. 
Due to the $\taumod$-regularity assumption on $(x_n)$, one finds (for all sufficiently large $n$) a unique face $\tau_n$ of type $\taumod$ in $\Flagt$, such that $\xi_n$ belongs to the open star $\ost(\tau_n)$ of $\tau_n$. By the definition, $x_n\to \tau$ (the sequence 
$(x_n)$ {\em flag-converges} to $\tau$)  if and only if  $\tau_n\to \tau$ in $\Flagt$. 

If $(x_n)$ is uniformly $\taumod$-regular (i.e., 
$\Theta$-regular for a compact $\Theta\subset \ost(\taumod)$) one can also describe flag-convergence $x_n\to\tau$ as follows. 
First, note that a diverging sequence $x_n\in X$ converges to $\xi\in \geo X$ with respect to the visual topology on $\ol{X}$ if and only if 
the sequence $(\xi_n)$ defined above converges to $\xi$ in the visual topology on $\geo X$. Of course, the sequence $(\xi_n)$ need not converge, but (by compactness of $\geo X$) it has convergent subsequences. In view of the $\Theta$-regularity of $(x_n)$,  all subsequential limits of $(\xi_n)$ in $\geo X$  (equivalently, of $(x_n)$ in $\ol{X}$) are $\Theta$-regular points in $\geo X$. Then $(x_n)$ flag-converges to $\tau\in \Flagt$ if and only if the accumulation set of $(x_n)$ in $\geo X$ is contained in $\ost(\tau)$.


\section{Regular and Anosov subgroups} \label{sec:regular}

{\bf Regular subgroups.} In what follows, we fix an $\iota$-invariant face $\taumod$ of $\simod$. 
(For the symmetric spaces appearing in this paper, the $\iota$-invariance condition is automatically satisfied since $\iota=\id$.) 
Importance of this invariance assumption comes from the fact that we will be interested in accumulation points in $\ol{X}^{\taumod}$ of $\Ga$-orbits of $\taumod$-regular subgroups $\Ga< G$. For a typical element $\gamma\in \Gamma$, if a sequence $(\gamma^n)_{n\in \N}$ is $\taumod$-regular, then  the inverse sequence $(\gamma^{-n})_{n\in \N}$ is $\iota\taumod$-regular. Hence, to have a satisfactory  theory, it makes sense to assume that $\taumod=\iota\taumod$.

\begin{rem}
We must also note that the notion equivalent to $\taumod$-regularity of subgroups $\Ga< G$ and the $\taumod$-limit set 
was first introduced by  Benoist in his highly influential work \cite[section 3.6]{Benoist}. For the benefit of an interested reader, 
his notation for the limit set was $\Lambda_\Gamma$.   
\end{rem}

We refer the reader to \cite{coco15} and \cite{manicures} for the detailed discussion of 
{\em $\taumod$-regular} discrete subgroups $\Ga< G$ and their $\taumod$-limit sets (denoted  $\Lat(\Ga)$ in our papers), 
which are certain closed $\Ga$-invariant subsets of $\Flagt$.

Below we review the notions of regularity and limit sets. A (necessarily discrete) subgroup $\Gamma< G$ is said to be {\em $\taumod$-regular} if every sequence of distinct elements $\gamma_n\in \Gamma$ is $\taumod$-regular. Similarly, one defines 
{\em uniformly $\taumod$-regular} subgroups of $G$. For instance, if $X$ has rank 1, then $\Delta$ is 1-dimensional, hence, uniform regularity of a subgroup is equivalent to discreteness. 

We next turn to the discussion of limit sets. Following \cite{Benoist}, for a discrete (not necessarily regular) subgroup $\Ga< G$ we define the {\em visual limit set} $\La(\Ga)\subset \geo X$ as the accumulation set of one (equivalently, every) $\Ga$-orbit $\Ga x\subset X$ 
with respect to the visual compactification of $X$. The next lemma is an immediate consequence of Lemma \ref{lem:horoball}: 

\begin{lemma}\label{lem:half}
Let $\Ga< G$ be a discrete subgroup preserving a horoball $Hbo\subset X$ centered at a point $\xi\in \geo X$. Then 
$$
\La(\Ga)\subset \bar{B}(\xi, \pihalf),
$$
the closed ball in $\tits X$, centered at $\xi$, of radius $\pihalf$ with respect to the Tits metric. 
\end{lemma}

\medskip
The $\taumod$-limit set $\Lat(\Ga)$ of a $\taumod$-regular subgroup $\Ga< G$ is the accumulation set in $\Flagt\subset \ol{X}^{\taumod}$ of some (equivalently, every) orbit $\Ga x\subset X$. In other words, 
$\tau\in \Lat(\Ga)$ if and only if there exists a sequence $(\ga_n)$ in $\Ga$ such that the sequence $(\ga_n(x))$ flag-converges to $\tau$. 
Since flag-convergence is independent of the base-point, $\Lat(\Gamma)$ is a closed $\Ga$-invariant subset of $\Flagt$. 

By the construction, since $\Lat=\Lat(\Ga)$ is $\Ga$-invariant, so is $\Th(\Lat)\subset \Flagt$ 
for every $\taumod$-invariant thickening $\Th=\Th_I$.  Since $\Lat$ is closed in $\Flagt$, so is $\Th(\Lat)$.  If $\Ga$ is uniformly $\taumod$-regular then $\Lat(\Ga)$ has an alternative description:
\begin{equation}\label{eq:lat} 
\Lat(\Ga)= \{\tau\in \Flagt : \ost(\tau)\cap \La(\Ga)\ne \emptyset\}, 
\end{equation}
cf. the alternative description of flag-convergence in the end of the previous section. 

\begin{cor}\label{cor:half}
Under the hypotheses of Lemma \ref{lem:half}, assume also that $G$ is a simple Lie group of type $B_2$ (hence, $\tits X$ is a graph with edges of length $\pi/4$), $\taumod$ is one of the two vertices of $\simod$, $\xi$ is a vertex of type $\taumod$, and $\Ga< G$ is a uniformly $\taumod$-regular subgroup. Then 
$$
\Lat(\Ga) \subset \bar{B}(\xi, \pihalf)\cap \Flagt\subset \tits X. 
$$
\end{cor}
\proof Note that if $\eta\in \tits X$ is a $\taumod$-regular point, $\tau\in \Flagt\subset \tits X$, then $\eta\in \ost(\tau)$  if and only if 
$\tangle(\eta, \tau)< \frac{\pi}{4}$. By Lemma \ref{lem:half}, $\La(\Ga)\subset \bar{B}(\xi, \pihalf)$. By combining these facts with \eqref{eq:lat}, we obtain  
$$
\Lat(\Ga) \subset \bigcup_{\eta\in \La(\Ga)} B(\eta, \frac{\pi}{4}) \cap \Flagt \subset B(\xi, \frac{3\pi}{4}) \cap \Flagt \subset 
\bar{B}(\xi, \pihalf)\cap \Flagt. \qed 
$$

\medskip 
A key result used in this paper is Theorem 6.13 from \cite{coco15}:

\begin{thm}\label{thm:disc}
Let $\Th$ be a fat thickening. Then for every $\taumod$-regular subgroup $\Ga< G$, the $\Ga$-action on 
$$
\Omega_{{\scriptsize \Th}}(\Ga):= \Flagt \setminus \Th(\Lat(\Ga))$$
is properly discontinuous. 
\end{thm}


{\bf Anosov subgroups.} 
An important class of $\taumod$-regular discrete subgroups $\Ga< G$ consists of {\em $\taumod$-Anosov subgroups}. Anosov representations $\Ga\to G$, whose images are Anosov 
subgroups, were first introduced in \cite{Labourie} for fundamental groups of closed manifolds of negative curvature, then in \cite{GW} for arbitrary hyperbolic groups; we refer the reader to our papers \cite{anolec, mlem, bordif}, for a simplification of the original definition as well as for alternative definitions and to \cite{manicures, anosov} for surveys of the results.

\medskip
Instead of a detailed discussion of  Anosov subgroups, we limit ourselves here to a brief description of their key properties used in this paper. Firstly, suppose that $H$ is a rank one Lie group and $X_H$ be the corresponding rank one symmetric space (the reader can assume that $H=O(n-1,1)$  and $X_H$ is the hyperbolic $n-1$-space $\H^{n-1}$). Then the Tits topology on  $\geo X_H$ is discrete. Accordingly, there is only one type of visual boundary simplices $\taumod=\taumod^H$ and, as we noted earlier, a subgroup $\Ga< H$ is discrete if and only if it is  $\taumod$-regular. The $\taumod$-limit set $\Lat(\Ga)\subset \geo X_H$ is the visual limit set $\La(\Ga)$. A subgroup $\Ga< H$ is Anosov (more precisely, $\taumod^H$-Anosov) if and only if it is {\em convex-cocompact}, equivalently, if it is discrete, finitely-generated and one, equivalently, every, orbit map $\Gamma\to X_H$ is a quasiisometric embedding of $\Gamma$ (equipped with a word-metric) to the symmetric space $X_H$. See for instance,  Theorem 1.1 in \cite{anolec} and also \cite{Bowditch2}. 

\medskip
Now consider the case of discrete subgroups of a semisimple Lie group $G$ without any restriction on rank; $X=G/K$ is the associated symmetric space. Suppose that $\taumod$ is an $\iota$-invariant face of $\simod$. Below are two of the many characterizations of $\taumod$-Anosov subgroups $\Ga< G$ given in \cite{anolec, mlem, bordif}:

\begin{thm}\label{thm:anosov}
The following are equivalent for a subgroup $\Ga< G$:

1. $\Ga$ is Gromov-hyperbolic, $\taumod$-regular (as a subgroup of $G$), any two distinct limit points in $\Lat(\Ga)\subset \Flagt$ 
are antipodal and there exists an equivariant homeomorphism $\beta: \geo \Ga\to \Lat(\Ga)$. Here $\geo \Ga$ is the Gromov-boundary of $\Gamma$.  The map $\beta$ is called {\em the boundary map} of $\Ga$. 

2. $\Ga$ is finitely generated, uniformly $\taumod$-regular (as a subgroup of $G$) and is {\em undistorted}, i.e. one (equivalently, every) orbit map $o_x: \Gamma\to \Gamma x\subset X$ is a quasiisometric embedding. 

3. $\Ga< G$ is $\taumod$-Anosov. 
\end{thm}

\medskip 
{\bf Images of rank 1 Anosov subgroups in higher rank Lie groups.}  Suppose that $G$ is a semisimple Lie group (the reader can assume that $G=O(n,2)$) and $H\to G$ is an embedding of Lie groups (the reader can think of the natural inclusion $O(n-1,1)\to O(n,2)$; it the one given by the composition of the embeddings $O(n-1,1)\to G_{L, \hat{L}}\to G$ discussed in the next section). For simplicity of the discussion (and because it is true in the main example of interest), we assume that the opposition involution $\iota$ of the group $G$ is the identity map. Let $X=G/K$ be the symmetric space of $G$, $X_H$ is the symmetric space of $H$ and let $X_H\to X$ be a totally-geodesic embedding equivariant with respect to the embedding $H\to G$.  (In the context of $H=O(n-1,1)< G=O(n,2)$, we will discus the embedding $X_H\to X$ in Section \ref{sec:main}.) The embedding $X_H\to X$ induces an isometric embedding of Tits boundaries $\tits X_H\to \tits X$ (this embedding is not in general simplicial, but it will be simplicial in the case of interest in this paper); we will identify $\geo X_H$ with its image in $\geo X$. Accordingly, for every point $\eta\in \tits X_H$, there exists a unique smallest simplex $\tau:=\xi(\eta)$ 
in $\tits X$ containing $\eta$. (In other words, $\eta$ is a generic point of $\tau$.) 
All the simplices $\tau=\xi(\eta)$ have the same type, which we denote $\taumod$. (In the case of interest, we will see that $\xi(\eta)$ is always a vertex of the type of an isotropic line, i.e. an element of the flag-manifold $\F_1$. Hence, in this case $\xi$ is the identity embedding.) The map $\xi: \geo X_H\to \Flagt$ is continuous, where $\geo X_H$ is equipped with the visual topology. 
It follows from the main definition of the $\taumod$-regularity and $\taumod$-limit set that for a discrete subgroup $\Ga< H$, its image in $G$ (also denoted $\Ga$) is uniformly $\taumod$-regular and that $\Lat(\Ga)=\xi(\La(\Ga))$, where $\La(\Ga)$, as we noted earlier, is the limit set of $\Ga$ in the visual boundary of $X_H$. 
Furthermore, it follows immediately from {\em every} characterization of $\taumod$-Anosov subgroups of $G$ 
given in \cite{anolec, mlem} (see for instance Theorem \ref{thm:anosov} above) that if $\Ga< H$ is convex-cocompact, then $\Ga< G$ is 
$\taumod$-Anosov. This fact was first observed by Labourie in \cite[Proposition 3.1]{Labourie} in the {\em Fuchsian case} and then in \cite[Proposition 4.7]{GW} in full generality.   We summarize these observations in the following proposition:

\begin{prop}\label{prop:embedding}
Let $G$ be a semisimple Lie group, $H< G$ is a rank 1 simple Lie subgroup, let $X_H\to X$ be a totally-geodesic embedding of the associated symmetric spaces, equivariant with respect to the embedding $H\to G$.  Then there exists a model face $\taumod$ of $\tits X$ such that the following hold for every discrete subgroup $\Ga< H$:

1. The image of $\Ga$ in $G$ is uniformly $\taumod$-regular. 

2. There exists a $\Ga$-equivariant homeomorphism $\beta: \La(\Ga)\to \Lat(\Ga)\subset \Flagt$ sending each $\la\in \La(\Ga)\subset \geo X_H$ to the unique simplex of type $\taumod$ in $\tits X$ containing $\la\in \geo X_H\subset \geo X$. 

3. If $\Ga< H$ is convex-cocompact, then $\Ga< G$ is $\taumod$-Anosov.  
\end{prop}

Note that the map $\beta$ here is the restriction of the map $\xi$ to $\La(\Ga)\subset \geo X_H$. It can be identified with the boundary map of the Anosov subgroup $\Ga< G$ as in Theorem \ref{thm:anosov} (the group $\Ga$ acts cocompactly the Gromov-hyperbolic space which is the closed convex hull $C$ of $\La(\Ga)$ in $X_H$ and, hence, $\geo \Ga$ can be identified with $\geo C=\La(\Ga)$).

\section{Lorentzian space $\R^{n-1,1}$ as an open Schubert cell in a partial flag-manifold of the group $G=O(n,2)$} \label{sec:identification}

In this section we will construct an equivariant identification of the Lorentzian space $\R^{n-1,1}$ 
with an open Schubert cell in a partial flag-manifold $\F_1$ of the group $G=O(n,2)$, namely, the space of isotropic lines in $V=\R^{n,2}$.

\medskip 
Consider the group $G=O(n,2)$ and its symmetric space $X=G/K$, $K=O(n)\times O(2)$. 
The group $G$ has two partial flag-manifolds: the Grassmannian $\F_1$ of isotropic lines and another partial 
flag  manifold $\F_2$ of isotropic planes in $V=\R^{n,2}$, where the quadratic form on $V$ is
$$
q= x_1 y_1 + x_2 y_2 + z_1^2+....+ z_n^2. 
$$
We will use the notation $\<\cdot, \cdot\>$ for the associated bilinear form on $V$. 

In the paper we will be using the  the incidence geometry interpretation of $\tits X$,   the Tits boundary  of the symmetric space 
of the group  $G=O(n,2)$.  The Tits boundary $\tits X$ (as a spherical building) has the structure of a metric bipartite graph whose vertices are labelled {\em lines} and {\em planes}, these are the elements of  
$\F_1$ and $\F_2$ respectively.  Two vertices $L\in \F_1$ and $p\in \F_2$ are connected by an edge iff the line $L$ is contained in the plane $p$. The edges of this bipartite graph have length $\pi/4$. We refer the reader to \cite{Brown},  \cite{Garrett} and \cite{Tits}. 

The group $G$ acts simply transitively on the set of edges of $\tits X$ and we can identify the quotient $\tits X/G$ with $\simod$, the model spherical chamber of $\tits X$. Thus $\simod$ is a circular segment of the length $\pi/4$. This segment has two vertices, one of which we denote $\taumod$, this is the one which is the projection of $\F_1$. The flag-manifold $\F_1$ is the quotient $G/P_L$, where $P_L$ is the stabilizer of an isotropic line $L$ in $G$; this flag-manifold is $n$-dimensional.

Recall that two vertices of $\tits X$ are opposite iff they are within Tits distance $\pi$ from each other. 
In terms of the incidence geometry of the vector space $(V,q)$, two lines $L, \hat L\in \F_1$ are opposite iff  
the restriction of $q$ to $\span(L, \hat L)$ 
is nondegenerate, necessarily of the type $(1,1)$. 
Two lines $L, L'\in \F_1$ are within Tits distance $\pi/2$ iff they span an isotropic plane in $V$. 

Consider a subgroup $P_L< G$; it is a maximal parabolic subgroup of $G$; 
let $U< P_L$ be the unipotent radical of $P_L$. Choosing a line  $\hat L$ opposite to $L$, defines a semidirect product decomposition $P_L= U\rtimes G_{L,\hat L}$, where $G_{L,\hat L}$ is the stabilizer in $P_L$ of the line $\hat{L}$; equivalently, it is the stabilizer of the {\em parallel set} $P(L, \hat L)$. \footnote{The parallel set 
$P(L, \hat L)$ splits isometrically as the product $l\times \H^{n-1}$, where $\H^{n-1}$ is the {\em cross-section} of $P(L, \hat L)$.} 
 This subgroup is the intersection
$$
G_{L,\hat L}= P_L\cap P_{\hat L}. 
$$

The orthogonal complement $V_{L,\hat L}\subset V$ of the anisotropic plane $\span(L, \hat L)$ 
is invariant under $G_{L,\hat L}$, hence, 
\begin{equation}\label{eq:splitting} 
G_{L,\hat L} \cong \R^\times \times O( V_{L,\hat L} , q|_{V_{L,\hat L}})\cong  \R^\times \times O(n-1,1). 
\end{equation}
The subgroup $\R_+< \R^\times$ acts via transvections along geodesics in the symmetric space $X$ connecting $L$ and $\hat L$. The group $G_{L,\hat L}$ acts on both $(V', q')= (V_{L,\hat L} , q|_{V_{L,\hat L}})$ and on $U$, where the action of $\R_+$ on   $V'=V_{L,\hat L}$ is trivial. In order to simplify the notation, we set
$$
O(q')= O(V', q'). 
$$
In terms of linear algebra, $\R_+< \R^\times$ is the identity component of the orthogonal group 
$$O(\span(L, \hat{L}), q|_{\span(L, \hat{L})})\cong O(1,1).$$

We will use the notation 
$$
G'_L:= U\rtimes O(q')< P_L.$$ 
This subgroup is the stabilizer in $P_L$ of horoballs in $X$ centered at $L$. 

\medskip
Our next goal is to describe Schubert cells in the Grassmannian $\F_1$. 
We fix $L\in \F_1$ and define the subvariety $Q_L\subset \F_1$  consisting 
of all (isotropic) lines $L'\subset V$ such that $\span(L, L')$ is isotropic (the line $L$ or an isotropic plane). 
 In terms of the Tits' distance, $Q_L \setminus \{L\}$ consists of lines $L'\in \F_1$ within distance $\pihalf$ from $L$. The complement 
$$
L^{opp}=\F_1 \setminus  Q_L$$
consists of lines opposite to $L$. The group $P_L$ acts transitively on $\{L\}$, $Q_L \setminus \{L\}$ and $L^{opp}$ and each of these subsets is an open Schubert cell of $\F_1$ with respect to $P_L$ and we obtain the $P_L$-invariant Schubert cell decomposition 
$$
\F_1 = \{L\} \sqcup (Q_L \setminus \{L\}) \sqcup L^{opp}. 
$$

We next describe $Q_L$ more geometrically. A vector $v\in V$ spans an isotropic subspace with $L$ iff 
$v\in L^\perp$ and satisfies the quadratic equation 
$q(v)=0$.  Since we are only interested in nonzero vectors $v\ne 0$ and their spans $\span(v)$, we obtain the natural identification
$$
Q_L \cong \P (q^{-1}(0) \cap L^\perp ), 
$$ 
the right hand-side is the projectivization a singular quadric hypersurface in $L^\perp$.  Thus, $Q_L$ is a (projective) singular quadric  
and $L\in Q_L$ is the unique singular point of the  $Q_L$. 


\medskip 
In the next lemma, by an {\em ellipsoid} in a real projective space $\R \P^{k-1}$ we mean the projectivization $E$ of a quadric in $\R^k$ given by a  quadratic form of signature $(k-1,1)$. (The reason for the name is that in a suitable affine patch in $\R \P^{k-1}$, $E$ becomes an ellipsoid.) 

\begin{lem}\label{lem:L4}
Given two opposite isotropic lines $L, \hat L$, the intersection of the quadrics  
$$
E= E_{L, \hat L}:= Q_L \cap Q_{\hat L}
$$
is an ellipsoid in  $\P  (L^\perp \cap \hat L^\perp)$. In particular, $E\cap \{L, \hat{L}\}=\emptyset$. 
\end{lem}
\proof As before, let $V'\subset V$ denote the codimension two subspace orthogonal to both $L, \hat L$. Then 
each $L'\in E$ is spanned by a vector $v\in V'$ satisfying the condition $q(v)=0$. In other words, $E$ is the 
projectivization of the quadric 
$$
\{v\in V': q(v)=0\}.
$$
i.e. is an ellipsoid, since $q$ restricted to $V'$ has signature $(n-1,1)$. \qed 

\medskip Our next goal is to (equivariantly) identify the open cell $L^{opp}$ with the $n$-dimensional Lorentzian affine space $\R^{n-1,1}$ (where a chosen $\hat L\in L^{opp}$ will serve as the origin), so that the group $P_L$ is identified with the group of Lorentzian similarities, where the simply-transitive action $U\acts L^{opp}$ is identified with the action of the full group of translations of $\R^{n-1,1}$. In particular, 
for now, and until the end of the proof of Proposition \ref{prop:P1}, we fix isotropic opposite lines $L$ and $\hat{L}$. After the end of the 
proof of Proposition \ref{prop:P1} we will allow the line $\hat{L}$ to vary.  

We fix nonzero vectors $e\in L$, $f\in \hat L$ such that $\<e, f\>=1$. Then 
$$
V= \span(e) \oplus \span(f) \oplus V'.
$$
We obtain an epimorphism $\eta: P_L\to O(q')$ by sending $g\in P_L$ first to the restriction $g|L^\perp$ and then to the projection of the latter to the quotient space $V'\cong L^\perp/L$ (the quotient of $L^\perp$ by the  null-subspace of $q|L^\perp$).  
Hence, the kernel of this epimorphism is precisely the solvable  radical $U\rtimes \R_+$ of $P_L$.

For each $v'\in V'$ we define the linear transformation (a shear) 
$s=s_{v'}\in GL(V)$ by its action on $e, f$ and $V'$:

\begin{enumerate}
\item $s(e)=e$. 

\item $s(f)= - \half q(v') e + f + v'$. 

\item For $w\in V'$, $s(w)= w - \< v', w\> e$. 
\end{enumerate}

The next two lemmata are proven by straightforward calculations which we omit: 

\begin{lem}
For each $s=s_{v'}$ the following hold:

1. $s\in P_L$.

2. $s$ lies in the kernel of the homomorphism $\eta: P_L\to GL(V')$ and is unipotent.  
In particular, $s\in U$ for each  $v'\in V'$. 
\end{lem}

 \begin{lem}
The map $\phi: v'\mapsto s_{v'}$ is a continuous monomorphism $V'\to U$, where we equip the vector space $V'$ with the additive group structure. 
\end{lem}

Since $U$ acts simply transitively on $L^{opp}$, it is connected and has dimension $n$. Therefore, the monomorphism $\phi$ is surjective and, hence, a continuous isomorphism. Thus, $\phi$ determines a homeomorphism $h: V'\to L^{opp}$ 
$$
h: v'\mapsto s_{v'}(\hat L)= span\left(- \half q(v') e + f + v'\right), 
$$
so that in particular
$$
h(0)= \hat L. 
$$

The group $G_{L,\hat{L}}\cong \R^\times \times O(V', q')$ acts on both $L^{opp}$ and on $U$ (via conjugation).  
The center of $G_{L,\hat{L}}$ acts on $V'$  trivially while its action on $U$ is via a nontrivial character.

\begin{prop}\label{prop:P1}
The map $h$ is equivariant with respect to these two actions of $O(V',q')$.  
\end{prop}
\proof Consider a linear transformation $A\in O(V',q')$; as before, we identify $O(V',q')$ with a subgroup of $O(V, q)$ fixing $e$ and $f$.   
For an arbitrary $v'\in V'$ we will verify that 
$$
s_{A v'}= A s_{v'} A^{-1}. 
$$
It suffices  to verify this identity on the vectors $e, f$ and arbitrary $w\in V'$. We have:  

1. For each $v'\in V'$, $s_{v'}(e)=e$, while $A(e)=A^{-1}(e)=e$. It follows that 
$$
e= s_{A v'}(e)= A s_{v'} A^{-1}(e)= e. 
$$ 

2. 
$$
s_{A v'}(f)= - \half q(A v') e + f + A v'= - \half q(v') e + f + Av'
$$
while (since $Ae=e$, $Af=f$) 
$$
A s_{v'} A^{-1}(f)= A s_{v'}(f)   = A (- \half q(v') e + f + v'   ) =  - \half q(v') e + f + Av'. 
$$

3. For $w\in V'$,
$$
s_{A v'} (w)= w- \< A v', w\> e= w -  \< v',  A^{-1} w\> e,
$$
while
$$
A s_{v'} A^{-1} w= A s_{v'}( A^{-1} w)= A ( A^{-1} w - \< v', A^{-1} w\> e)= w- \< v', A^{-1} w\> e.  \qed 
$$

In view of this proposition we will  identify $V'$ with the open Schubert cell $L^{opp}$, which, in turn, enables us to 
use Lorentzian geometry to analyze $L^{opp}$ and, conversely, to study discrete subgroups of $P_L$ using 
Theorem \ref{thm:disc} 
on domains of discontinuity of $\taumod$-regular group actions on the flag-manifold $\F_1$.  
Under the identification $V'\cong L^{opp}$, for each 
${L}'\in L^{opp}$ (in particular, for $L'=\hat{L}$), the conic $Q_{L'}\cap L^{opp}$ becomes a translate of the null-cone of the form $q'$ on $V'$ (see Lemma \ref{lem:L0} below) and the flag-manifold $\F_1$ becomes a compactification of $V'$ obtained by adding to it the ``quadric at infinity'' $Q_L$.

\begin{lem}\label{lem:L0}
For all $v'\in V'$, $q'(v')=0$ iff $q$ vanishes on $\span(f, h(v'))$, i.e. iff $h(v')\in Q_{\hat L}$.  In other words, 
$Q_{\hat L}\cap L^{opp}$ is the image under $h$ of the null-cone of $q'$ in the vector space $V'$. 
\end{lem}
\proof Since $f$ and $s_{v'}(f)$ (spanning the line $h(v')$) 
are null-vectors of $q$, the vanishing of $q$ on $\span(f, h(v'))$ is equivalent to the vanishing of 
$$
\< f, s_{v'}(f) \>= - \half q(v'). \qed 
$$

\begin{lem}\label{lem:L2}
For each neighborhood $N$ of $L$ in $Q_L$ there exists $L'\in L^{opp}$ such that $E_{L, L'}\subset N$.  
\end{lem}
\proof We pick $L_\infty\in \F_1$ opposite to $L$ and, as above, identify $L_\infty^{opp}$ with $(V',q')$. 
Then for a sequence $L_i\in L^{opp}_{\infty}$ contained in the, say, future light cone of $Q_L \cap L_{\infty}^{opp}$ 
and converging radially to $L$, the intersections of null-cones $E_{L,L_i}=Q_{L_i}\cap Q_L$ converge to $L$. Since $L_i\notin Q_L$, they are all opposite to $L$. Taking $L'=L_i$ for a sufficiently large $i$ concludes the proof. \qed

\medskip
For each subset $C\subset \F_1$, we define the {\em thickening} of $C$:
$$
\Th(C)= \bigcup_{L\in C} Q_L. 
$$
This notion of thickening is a special case of the one discussed in Section \ref{sec:prelim}:  If we restrict to a single apartment $a$  in the Tits building of $G$, then for the vertex $L\in a$, $\Th(L)\cap a= Q_L\cap a$ consists of three vertices within Tits distance   $\pihalf$ from $L$. Thus,  the thickening $\Th$ is maximal and, hence, {\em fat} (see Section \ref{sec:prelim}). 

\begin{lemma}\label{lem:L}
For every compact subset $C\subset \F_1$, the thickening $\Th(C)\subset \F_1$ is compact. 
\end{lemma}
\proof Compactness of thickenings (of closed subsets of general flag-manifolds)  is a general fact observed in \cite[p. 193]{coco15}, a proof can be found in \cite[Lemma 8.18]{bordif}, we add a proof here for the sake of completeness (it is the same as in \cite{bordif}). Compactness of $\Th(L)=Q_L$ for each $L\in \F_1$ is clear. Observe that  for $g\in G$, $g \Th(L)=\Th(gL)$. Consider a closed subset $C\subset \F_1$. Take a sequence $L_k\in C$ converging to $L_0\in C$. There exists a sequence $g_k\in K$ such that $g_k(L_1)=L_k$ for all $k\in \N$ (since the maximal compact subgroup $K< G$ acts transitively on $\F_1$). In view of compactness of the subgroup $K< G$, without loss of generality, we may assume that the sequence $g_k$ converges to some $g_0\in K$. Thus, the sequence of subsets $g_k(L_1)\subset \Th(C)$ converges to $g_0(L_1)$ with respect to the Hausdorff metric 
on the set of nonempty closed subsets of  $\F_1$.  At the same time, the sequence $g_k(L_1)=L_k$ converges to $L_0$, which implies that $L_0=g_0(L_1)$. Thus, the limit of the sequence of thickenings $g_k(\Th(L_1))$ equals the thickening $\Th(L_0)\subset \Th(C)$. It 
follows that $\Th(C)\subset \F_1$ is closed; compactness of $\F_1$ implies compactness of $\Th(C)$.

\begin{lem}\label{lem:L1}
For any two opposite isotropic lines $L, \hat{L}\in \F_1$ and each compact subset $C\subset Q_{\hat L}\cap L^{opp}$, 
the intersection $\Th(C)\cap L^{opp}$ is a proper subset of $L^{opp}$.  
\end{lem}
\proof Let $H\subset L^{opp}\cong V'$ be an affine hyperplane in $V'$ intersecting $Q_{\hat L}$ only at $\hat L$. Then 
$$
C':=\{L'\in H: Q_{L'} \cap C\ne \emptyset\}
$$
is compact in $H$. Next, observe that for $L_1, L_2\in \F_1$, $L_1\in Q_{L_2}\iff L_2\in Q_{L_1}$. Thus,  
every $L'\in H \setminus C'$ does not belong to $\Th(C)$. \qed

\begin{lem}\label{lem:L3}
For each compact $C\subset Q_L \setminus  \{L\}$ the thickening $\Th(C)$ is a proper compact subset of $\F_1$. 
\end{lem}
\proof Lemma \ref{lem:L2} implies that there exists $L_\infty\in L^{opp}$ such that $E_{L, L_{\infty}}$ is disjoint from $C$. 
Thus, $C$ is contained in $L^{opp}_{\infty}$. In particular, 
Lemma \ref{lem:L1} implies that $\Th(C)$ is a proper subset of $\F_1$. Compactness of $\Th(C)$ was proven in Lemma \ref{lem:L}. 
\qed

\section{Proof of the main theorem}\label{sec:main}

We continue with the notation introduced in the previous section. Consider the subgroups $G'_L < P_L < G$ with $G=O(n,2)$. 
The subgroup $H=O(n-1,1)< G$ stabilizes two opposite points in 
$\geo X$, which are isotropic lines $L, \hat{L}$ and, hence, preserves the parallel set $P(L, \hat{L})$ consisting of geodesics in $X$ asymptotic to both $L, \hat{L}$. This parallel set splits as the product ${\mathbb R}\times Y$, where $Y$ is a totally-geodesic symmetric subspace in $X$ (necessarily of rank one) and for each $y\in Y$ the product $\R\times \{y\}$ is one of the geodesics in $X$ asymptotic to    
$L, \hat{L}$. The subgroup $H$ preserves $P(L, \hat{L})$; it also necessarily preserves the product decomposition. The identity  
component  $H_0< H$ also necessarily preserves each $\{t\}\times Y$ (for otherwise, we obtain a nontrivial isometric action of $H$ on the real line). Moreover, since $H$ preserves both $L$ and $\hat{L}$, the entire group $H$ preserves each $\{t\}\times Y$. 
Pick a point $y\in Y$ and take a visual boundary point $\eta\in \geo Y$. We have two geodesic rays in $P(L, \hat{L})$ emanating from $y$: One is asymptotic to $L$, another (contained in $Y$) asymptotic to $\eta$. These rays are obviously contained in a 2-dimensional flat in $X$ and are orthogonal to each other. Hence, the Tits angle between $\eta, L$ equals $\pi/2$. Since the Tits boundary of $X$ is a bipartite graph with edge-length $\pi/4$, and $L$ is a vertex of this graph, the point $\eta$ is also a vertex and has the same type as $L$, 
i.e. the type of an isotropic line. Similarly, $\eta$ has the Tits distance $\pi/2$ from $\hat{L}$. 
Since the subgroup $H=O(n-1,1)< G=O(n,2)$ preserves each $\{t\}\times Y$, there exists a totally-geodesic isometric embedding of the symmetric space $X_H$ of $H$ into $\{t\}\times Y$. (Actually, \eqref{eq:splitting} implies that $X_H$ is the entire $\{t\}\times Y$ but we will not need this fact.) 

From now on, $\taumod$ is a vertex of the Tits building $\geo X$ which has the type of an isotropic line. 
In view of Proposition \ref{prop:embedding}, we conclude:

\begin{lem}
Let $\Ga< H$ be a discrete subgroup. Then the image of $\Gamma$ under the embedding $H\to G$ is $\taumod$-regular. Every $\taumod$-limit point $\eta$ of $\Gamma$ is a vertex of $\geo X$ of the type of an isotropic line, which is at the Tits distance $\pi/2$ from both $L, \hat{L}$. Accordingly, $\eta$ belongs to the intersection $Q_L\cap Q_{\hat{L}}$. 
\end{lem}

\begin{cor}\label{cor:Anosov}
If $\Gamma< H$ is a convex-cocompact subgroup, then its image in $G$ is $\taumod$-Anosov and its $\taumod$-limit set $\Lat(\Ga)$ 
is contained in $Q_L\cap Q_{\hat{L}}$. 
\end{cor}

We next consider the slightly more general case of uniformly $\taumod$-regular discrete subgroups $\Ga< P_L$: 

\begin{lem}\label{lem:L16}
The $\taumod$-limit set $\Lat(\Ga)$ of every uniformly $\taumod$-regular subgroup $\Ga< P_L< G$ is 
contained in $Q_L$.  
\end{lem}
\proof According to Corollary \ref{cor:half}, $\Lat(\Ga)\subset \bar{B}(L, \pihalf)\cap \F_1$. The latter intersection is $Q_L$ since both consist of isotropic lines $L'\subset V$ such that $\span(L, L')$ is an isotropic subspace of $V$. \qed 

\begin{prop}\label{prop:pd}
Suppose that $\Ga< G'_L$ is a $\taumod$-regular discrete subgroup whose $\taumod$-limit set does not contain $L$. Then 
$\Th(\Lat(\Ga))$ is closed in $\F_1$, $\Th(\Lat(\Ga))\ne \F_1$,  and the action 
$$
\Ga \acts \Om_{{\scriptsize \Th}}(\Ga)= \F_1 \setminus \Th(\Lat(\Ga)) 
$$
is properly discontinuous. 
\end{prop}
\proof Since $\Lat(\Ga)$ is a compact subset of $Q_L$, the first statement of the proposition is a special case of 
Lemma \ref{lem:L3}. The proper discontinuity statement is a special case of Theorem \ref{thm:disc}  
since the thickening $\Th$ is fat. \qed

\medskip 
We now describe certain  conditions on $\taumod$-regular discrete subgroups $\Ga< G'_L$ 
which will ensure that $\Lat(\Ga)$ does not contain the point $L$. Each subgroup $\Ga< G'_L$ has the {\em linear part} 
$\Ga_0$, i.e. its projection to $O(q')\cong O(n-1,1)$, which is identified with the semisimple factor of the stabilizer in $P_L$ of some $\hat L\in L^{opp}$. We now assume that:

\begin{itemize}
\item $\Ga_0$ is a convex-cocompact subgroup of $O(n-1,1)$.
\item The projection
$$
\ell: \Ga \to \Ga_0
$$ 
is an isomorphism. 
\end{itemize}

As we proved in Corollary \ref{cor:Anosov}, $\Ga_0< G$ is a ${\taumod}$-Anosov subgroup of $G$ and $\Lat(\Ga_0)\subset Q_L\cap Q_{\hat{L}}$. 
In particular, $\Lat(\Ga_0)$ does not contain $L$ by Lemma \ref{lem:L4}. 

Given a subgroup $\Ga_0< O(q')$, the inverse 
$\rho: \Ga_0\to \Ga$ to $\ell: \Ga\to \Ga_0$ is determined 
by a cocycle $c\in Z^1(\Ga_0, V')$ which describes the translational parts of the elements of $\Ga$:
$$
\rho(\ga): v\mapsto \ga v + c(\gamma), v\in V'\cong \R^{n-1,1}. 
$$
Pick some $t\in \R_+$; then $t c$ is again a cocycle corresponding to the conjugate representation $\rho^t$, where we identity $t\in \R_+$ with a central element of $G_{L,\hat L}$. Sending $t\to 0$ we obtain:
$$
\lim_{t\to 0} \rho^t = id,
$$ 
the identity embedding $\Ga_0\to O(n-1,1)< P_L$. In view of stability of Anosov representations (see \cite[Theorem 5.13]{GW} and \cite[Theorems 1.10, 1.11]{morse}, \cite[Corollary 6.14]{KKL}) we conclude that all representations $\rho^t$ are ${\taumod}$-Anosov and the $\taumod$-limit sets of $\Ga_t=\rho^t(\Ga_0)$ vary continuously with $t$; moreover,
$$
t \Lat(\Ga_{t_1})= \Lat(\Ga_{t_2})
$$
where $t=t_2/t_1$. In particular, 
$$
\Lat(\Ga) \subset Q_L \setminus \{L\}
$$
is a compact subset. Proposition \ref{prop:pd} now implies:

\begin{cor}
For each $\Ga$ as above, 
$$
Th(\Lat(\Ga))\ne \F_1 
$$
and the action
$$
\Ga \acts \Om_{{\scriptsize \Th}}(\Ga)= \F_1 \setminus Th(\Lat(\Ga))
$$
is properly discontinuous. 
\end{cor}

Thus, we proved that each discrete subgroup $\Ga< P_L$ as  above has nonempty domain of discontinuity in the vector space $V'$. Theorem \ref{thm:main} follows. \qed

\smallskip
{\bf Acknowledgements.} The first author was partly supported by the NSF grant  DMS-16-04241, 
by a Simons Foundation Fellowship, grant number 391602,  by Max Plank Institute for Mathematics in Bonn, 
as well as by KIAS (the Korea Institute for Advanced Study) through the KIAS scholar program. 
Much of this work was done during our stay at KIAS and we are thankful to KIAS for its hospitality. We are grateful to the referee of the paper for reading the paper carefully and for useful remarks and suggestions.

\noindent M.K.: Department of Mathematics, 
University of California, Davis, 
CA 95616, USA\\
email: kapovich@math.ucdavis.edu

and

\medskip 
\noindent Korea Institute for Advanced Study,\\ 
207-43 Cheongnyangri-dong, Dongdaemun-gu,\\ 
Seoul, South Korea\\

\noindent B.L.: Mathematisches Institut,
Universit\"at M\"unchen, 
Theresienstr. 39, 
D-80333 M\"unchen, Germany, 
email: b.l@lmu.de

\end{document}